\documentclass[10pt]{article}
\usepackage{amssymb,amsmath,makeidx,amsfonts,amsthm,verbatim,mathrsfs,euscript,eufrak}
\input xy
\xyoption{all}
\newcommand{\R}{I\!\!R}

\newcommand{\h}{I\!\!H}

\newcommand{\C}{I\!\!\!\!C}

\newcommand{\be}{\begin{enumerate}}
\newcommand{\ee}{\end{enumerate}}
\newcommand{\bq}{\begin{eqnarray*}}
\newcommand{\eq}{\end{eqnarray*}}

\begin{document}
\newcommand{\disp}{\displaystyle} \baselineskip 18pt
\thispagestyle{empty}
\begin{center}
{\bf A NEW SPHERICAL HARMONICS ON THE HEISENBERG GROUP}
\ \\
\ \\
{\bf M.E. Egwe\\
Department of Mathematics,\\
Universty of Ibadan\\
$murphy.egwe@mail.ui.edu.ng$}\\
\ \\
\end{center}
{\bf Abstract:}
\emph{Let $\h_n$ be the $(2n+1)$-dimensional  Heisenberg group. and let ${\cal L}_\alpha$ be the sublaplacian of the Lie algebra of $\h_n$ A new spherical harmonics with its orthogonal polynomial properties is presented for the group.} 

\section{Preliminaries}

The Heisenberg group (of order $n$), $\h_{n}$ is a
noncommutative nilpotent Lie group whose underlying manifold is $\C^{n}\times \R$
with coordinates $(z,t)=(z_{1},z_{2},...,z_{n},t)$ and group law given by
$$ (z,t)(z',t')=(z+z',t+t'+2\hbox{Im}\;( z.z'),\;\; \mbox{where}\;\;
 z.z'=\sum_{j=1}^{n}z_{j}\bar z_{j}\;\;\;\;z\in\C^{n},\;\;t\in\R. $$
Setting $z_{j}=x_{j}+y_{j}$, then $(x_{1},x_{2},...,x_{n},y_1,y_2,\cdots,y_n,t)$ forms a real coordinate
system for $\h_{n}$.
In this coordinate system, we define the following vector fields:
$$ X_{j}=\frac{\partial}{\partial x_{j}}+2y_{j}\frac{\partial}{\partial t},\;\; Y_{j}=\frac
{\partial}{\partial y_{j}}-2x_{j}\frac{\partial}{\partial t},\;\;T=\frac{\partial}{\partial t}. $$ It is clear from [1] that $\{X_{1},X_{2},...,X_{n},Y_{1},Y_{2},...,Y_{n},T\}$ is a basis for the left invariant vector fields on $\h_{n}.$ These vector fields span the Lie algebra $\mathfrak{h}_n$ of $\h_n$ and
the following commutation relations hold:
 $$[Y_{j},X_{k}]=4\delta_{jk}T,\;\;\; [Y_{j},Y_{k}]=[X_{j},T]=[Y_{j},T]=0.$$
Similarly, we obtain the complex vector fields by setting
$$\left.\begin{array}{lr}
Z_j =\frac{1}{2}(X_j-iY_j) = \frac{\partial}{\partial z_j} +
i\bar{z}\frac{\partial}{\partial t}\\
\ \\
\bar{Z}_j =\frac{1}{2}(X_j+iY_j) = \frac{\partial}{\partial
\bar{z}_j} -
iz\frac{\partial}{\partial t}\end{array}\right\}.\eqno(1.1)$$
In the complex coordinate, we also have the commutation relations
$$[Z_j,\bar{Z}_k] = -2\delta_{jk}T,\;\;[Z_j,Z_k] = [\bar{Z}_j,\bar{Z}_k]= [Z_j,T]=[\bar{Z},T] = 0.$$
The Haar measure on $\h_n$ is the Lebesgue measure
$dzd\bar{z}dt$ on $\C^n\times\R$  [2]. In particular, for $n=1$, we obtain the $3$-dimensional Heisenberg group $\h_1\cong\R^3$ (since $\C^n\cong\R^{2n})$. Hence $\h_n$ may also be referred to as $(2n+1)$-dimensional Heisenberg group, $[10],[11],[12]$.

One significant structure that accompanies the Heisenberg group is the family of dilations
$$\delta_{\pm \lambda}(z,t)=(\pm \lambda z,\pm \lambda^2t),\;\; \lambda >0 $$
This family is an automorphism of $\h_n.$ Now, if $\sigma:\C\rightarrow \C$ is an automorphism, there exists an induced automorphism, $\widetilde{\sigma} \in Aut(\h_n),$ such that
$$\widetilde{\sigma}(z,t)=(\sigma z, t).$$ For simplicity, assume that $\widetilde{\sigma}$ and $\sigma$ coincide. Thus we may simply assume that if $\sigma\in Aut(\h_n),$ we have $\sigma (z,t)=(\sigma z, t).$

An operator that occurs as an analogue (for the Heisenberg group) of the Laplacian\\
$\Delta=\disp\sum_{j=1}^n\frac{\partial^2}{\partial {x^2_j}}$ on $\R^n$ is denoted by ${\cal L}_\alpha$
where $\alpha$ is a parameter and defined by
$${\cal L}_\alpha=-\frac{1}{2}\sum_{j=1}^n(Z_j\bar{Z}_j+\bar{Z}_jZ_j)+i\alpha T,$$
where $\bar {Z}_j\;\mbox{and}\; Z_j$ are as defined in (1)
so that ${\cal L}_\alpha$ can be written as
$${\cal L}_\alpha=\frac{1}{4}\sum_{j=1}^n(X_{j}^2+Y_{j}^2)+i\alpha T.\eqno (1.2)$$
${\cal L}_\alpha$ is called the \emph{sublaplacian.}
${\cal L}_\alpha$ satisfies symmetry properties analogous to those of $\Delta$ on $\R^n.$ Indeed, we have that ${\cal L}_\alpha$
\be
\item[(i)] is left-invariant on $\h_n$
\item[(ii)] has degree $2$ with respect to the dilation automorphism of $\h_n$ and
\item[(iii)] is invariant under unitary rotations. \ee
Several methods for the determination of solutions, fundamental solutions of (2) and conditions for local solvability are well known [3,4,5].

The Heisenberg-Laplacian is a
subelliptic differential operator defined for $\alpha=0$ as
$\Delta_{\h_n}$ on $\h_n$ and denoted by ${\cal L}.$ It is obtained from the usual vector fields as\\
\begin{eqnarray*}
{\cal L}:=\Delta_{\h_n}&:=& \sum^n_{j=1}X_j\circ X_j+Y_j\circ Y_j\\
&=&\sum^n_{j=1}\frac{\partial^2}{\partial
x_j^2}+\frac{\partial^2}{\partial
y_j^2}+4y_j\frac{\partial^2}{\partial x_j\partial
t}-4x_j\frac{\partial^2}{\partial y_j\partial
t}+4(x_j^2+y_j^2)\frac{\partial^2}{\partial t^2}.\;\;\;\;\;\;\;\;\;\;\;\;\;\;\;\;\;\;\;\;\;\;\;\;\;\;
\hfill (1.3)\end{eqnarray*}
\section*{2. Spherical Harmonics}
On $\R^n$ following [1], [6], and [7].
The solutions of the Laplace's equation
$$\sum^n_{i=1}\frac{\partial^2u}{\partial x^2_i}=0\eqno (2.1)$$
are called harmonic function.  We are particularly interested in polynomial solutions.  Since any polynomial in the variables $x_1,x_2,\cdots,x_n$ is a sum of a finite number of homogeneous polynomials of different degrees, we concentrate on these.\\
\\
\textbf{Definition 2.1.1}~
A polynomial $H_m(x)$, which is homogeneous of degree $m$ and satisfies$(2.1)$ is called a harmonic polynomial.

Let $H_k(x)$ and $H_j(x)$ be homogeneous harmonic polynomials in $n$ variables of degree $k$ and $j$ respectively with $j\neq k$.  By Green's theorem,
$$\begin{array}{rcl}
0&=&\disp\int_{_{|x|\leq 1}}[H_j(x)\Delta H_k(x)-H_k(x)\Delta H_j(x)]dV\\
\ \\
&=&\disp\int_{_{|\xi|\leq 1}}\left[H_j(\xi)\frac{\partial}{\partial r}H(r\xi)|_{r=1}-H_k(\xi)\frac{\partial}{\partial r}H_j(r\xi)|_{_{r=1}}\right]dw(\xi)\end{array}\eqno(2.2)$$
where $\xi=x/|x|,\;|x|=r$ and $dw(\xi)$ is the invariant measure on the surface of the sphere (i.e., by using the fact that the normal derivative on the sphere is in the radial direction).  The homogeneity of $H_k(x)$ gives
$$\frac{\partial}{\partial r}H_k(r\xi)|_{r=1}=\frac{\partial}{\partial r}r^kH_k(x)|_{r=1}=kH_k(x).
\eqno (2.3)$$
Substituting $(2.3)$ in $(2.2),$ we obtain
$$(k-j)\int_{|\xi|=1}H_j(\xi)H_k(\xi)dw(\xi)=0.\eqno (2.4)$$\\
\\
\textbf{Definition 2.1.2}~  The functions $H_k(\xi)$, which are restrictions of homogeneous harmonic polynomials to the surface of the sphere in $\R^n$, are called spherical harmonics.  Sometimes, they are also called surface spherical harmonics and the $H_k(x)$ are called solid spherical harmonics.\\
A surface spherical harmonic of degree $n$ multiplied by the factor $r^n$ is called a solid spherical harmonic of degree $n$.

On the polar coordinates, spherical harmonics are denoted by
$$H_k(\xi)=Y_k(\theta,\phi)\eqno (2.5)$$
where $\theta=(\theta_1,\theta_2,\cdots,\theta_{n-2}).$\\

Let $\mathbb{C}(\R^n))$ be the space of real-valued continuous functions on the sphere\\ $|\xi|^2=\xi^2_1+\cdots+\xi^2=1.$  An inner product on this space can be defined by
$$\langle f,g\rangle = \int_{|\xi|=1}f(\xi)g(\xi)dw(\xi),\;\;f,g\in {\cal C}(\R^n).\eqno (2.6)$$
(If the functions are complex-valued, the conjugate of $g$ is used.)\\

We now give a version of spherical harmonics on the Heisenberg group following [52] in what follows.\\
\textbf{Definition 2.1.3}~The Heisenberg unit ball is defined by
$$\mathbb{B}_{\h_1}:=\{(z,t)\in \h_n:(|z|^4+t^2)^{\frac{1}{4}}<1\}$$
and is denoted $\mathbb{B}_{\h}(1).$\\
\\
\textbf{Definition 2.1.4}~ A solid ${\cal L}$-spherical harmonic of degree $m,\;\;m=0,1,2,\cdots$ is a polynomial in $z,\bar{z}$ and $t$, which is harmonic with respect to ${\cal L}$ (that is a function $u$ such that $Du=0$, where $D$ is an invariant differential on $\h_n$) and which is homogeneous of degree $m$ with respects to the Heisenberg dilation.

In what follows, we present an explicit computation of the homogeneous\\ ${\cal L}_0$-harmonic polynomial on $\h_1.$  That is, the solid ${\cal L}-$spherical harmonics and their boundary value on the boundary of the unit Heisenberg ball, $b\mathbb{B}_{\h}(1)$, the ${\cal L}_0$-spherical harmonics.\\
Mimicking  $\R^3$, we set
$$\begin{array}{rcl}
\rho&=&(|z|^4+t^2)^{\frac{1}{4}},\\
\\
x&=&\rho\sin^{\frac{1}{2}}\phi\cos\theta,\;\;y=\rho\sin^{\frac{1}{2}}\phi\sin\theta\\
\ \\
t&=&\rho^2\cos\phi,\end{array}$$
with $0<\rho<\infty,\;\; 0\leq\phi\leq\pi$ and $0\leq \theta\leq 2\pi$.
In other words,
$$\begin{array}{rcl}
z=|z|e^{i\theta}&=&\rho\sin^{\frac{1}{2}}\phi e^{i\theta}\\
\ \\
t+i|z|^2&=&\rho^2e^{i\phi}.\hspace{4.0in}(2.7)\end{array}$$
An explicit calculation yields
$$\begin{array}{rcl}
Z&=&\disp\frac{\partial}{\partial z}+i\bar{z}\frac{\partial}{\partial t}\\
\ \\
&=&e^{i\theta}\disp\left(\frac{1}{2}i\sin^{\frac{1}{2}}\phi e^{-i\phi}\frac{\partial}{\partial\rho}+\frac{\sin^{\frac{1}{2}}\phi}{\rho}\;e^{i\phi}\;
\frac{\partial}{\partial\phi}-\frac{1}{2}i\frac{1}{\rho\sin^{\frac{1}{2}}\phi}\frac{\partial}{\partial
\theta}\right).\end{array}$$
We shall also need
$$\frac{\partial}{\partial t}=\frac{\cos\phi}{2\rho}\frac{\partial}{\partial\rho}-\frac{\sin\phi}{\rho^ 2}\;\frac{\partial}{\partial\phi}.$$\\
In [8] , Folland and Stein introduced the operators ${\cal L}_{\alpha},\;\;\alpha\in\C$ given by
$${\cal L}_{\alpha}=\frac{1}{2}(Z\bar{Z}+\bar{Z}Z)+i\alpha\frac{\partial}{\partial t}$$
where
$$[Z,\bar{Z}]=Z\bar{Z}-\bar{Z}Z.$$
These operators behave like ${\cal L}_0.$ Thus, ${\cal L}_{\alpha}$ is left-invariant on $\h_1$ and homogeneous of degree two in the form
$${\cal L}_{\alpha}(f(RZ, R^2t))=R^2({\cal L}_{\alpha}f)(RZ,R^2t),\;\;R>0.$$
Furthermore, ${\cal L}_{\alpha}$ is hypoelliptic, i.e. ${\cal L}_{\alpha}u=f\in C^{\infty}\Longrightarrow u\in C^{\infty}$, and solvable if and only if $\alpha\neq \pm1,\pm 3,\pm 5,\cdots$.  When $\pm\alpha=1,2,3,5,\cdots,{\cal L}_{\alpha}$ has a fundamental solution given by $\Phi_{\alpha}(z,z',t-t'-\frac{1}{2}Im\;z.\bar{z}')=\delta_{z',t'}$ and
$$\Phi_{\alpha}(z,t)=C_{\alpha}(|z|^4+t^2)^{-\frac{1}{2}}\left(\frac{|z|^2+it}{|z|^2-it}\right)^{\alpha/2}$$
with
$$C_{\alpha}=\frac{1}{\pi^2}\Gamma\left(\frac{1+\alpha}{2}\right)\Gamma\left(\frac{1-\alpha}{2}\right).$$
When $\alpha=\pm1, \pm3,\pm5,\cdots$ one obtains the relative fundamental solutions as in [5, 10]\\
\\
\subsection*{2.2 The Spherical Harmonics for ${\cal L}_{\alpha}$}, $\pm\alpha= 1, 3, 5,\cdots$\\
We now compute the spherical harmonics for ${\cal L}_{\alpha}$.  We shall only assume the case where ${\cal L}_\alpha$ is soluble. i.e., $\pm\alpha=1,3,5,\cdots$.\\
\\
\textbf{Remark 2.2.1}~
In finding the homogeneous harmonic polynomials for ${\cal L}_{\alpha},$ notice that ${\cal L}_\alpha$
does not separate points in $(\rho,\theta,\phi)$ coordinates.  Therefore, to get an idea of what solid spherical harmonics for ${\cal L}_{\alpha}$ look like we first compute some explicitly for ${\cal L}_0$ as shown below.

Let ${\cal H}^{(0)}_m$ denote the set of homogeneous harmonic polynomials of degree $m$ for ${\cal L}_0$.  Then
$$\begin{array}{rcl}
{\cal H}^{(0)}_0&=&\{1\},\\
\ \\
{\cal H}^{(0)}_1&=&\{z,\bar{z}\},\\
\ \\
{\cal H}^{(0)}_2&=&\{z^2,\bar{z}^2,t\},\\
\ \\
{\cal H}^{(0)}_3&=&\left\{z^3,\bar{z}^3,z(|z|^2-2it),\bar{z}(|z|^2+2it)\right\},\\
\ \\
{\cal H}^{(0)}_4&=&\left\{z^4,\bar{z}^4,z^2(|z|^2-\frac{3}{2}it),\bar{z}^2(|z|^2+\frac{3}{2}it), |z|^4-2t^2\right\},\\
\vdots&&\mbox{e.t.c.}\end{array}$$ In general, we have
$${\cal H}^{(0)}_n=\{z^n,\bar{z}^n,z^{n-2}(|z|^2-\frac{n-1}{2}it),
\bar{z}^{n-2}(|z|^2+\frac{n-1}{2}it),(z^n-\frac{n}{2}it^2),...\}$$
(considering the homogeneity with respect to dilations).

Introducing spherical coordinates, we obtain
$$\begin{array}{rcl}
{\cal H}^{(0)}_0&=&\{1\},\\
\ \\
{\cal H}^{(0)}_1&=&\left\{e^{i\theta}\rho\sin^{1/2}\phi,\;e^{-i\theta}\rho\sin^{1/2}\phi\right\},\\
\ \\
{\cal H}^{(0)}_2&=&\left\{e^{2i\theta}\rho^2\sin\phi,\;e^{-2i\theta}\rho^2\sin\phi,\rho^2\cos\phi\right\},\\
\ \\
{\cal H}^{(0)}_3&=&\left\{e^{3i\theta}\rho^3\sin^{3/2}\phi,\;e^{-3i\theta}\rho^3\sin^{3/2}\phi,e^{i\theta}\rho^3
\sin^{3/2}\phi\left(1-2i\frac{t}{|z|^2}\right)\right\},\\
\ \\
{\cal H}^{(0)}_4&=&\disp\left\{e^{4i\theta}\sin^2\phi,\;e^{-4i\theta}\rho^4\sin^2\phi,e^{2i\theta}\rho^4
\sin^2\phi\left(1-\frac{3}{2}i\frac{t}{|z|^2}\right)\right.,\\
\ \\
&&e^{-2i\theta}\rho^4\sin^2\phi,\left(1+\frac{3}{2}i\frac{t}{|z|^2}\right),\rho^4
\sin^2\phi\left(1-2\left(\frac{t}{|z|^2}\right)^2\right\},\\
\vdots&&\mbox{e.t.c.}\end{array}$$
From these formulas, we can make a guess at the form of the general solid spherical harmonic of ${\cal L}_0$, and, also, of that of ${\cal L}_{\alpha}$ in the following result.\\
\\
\textbf{Proposition 2.2.2}~
All solid spherical harmonic of ${\cal L}_{\alpha}$ of degree $m,$\\ $m = 0,1,2,\cdots$ are sum of $\h_1$-homogeneous ${\cal L}_{\alpha}$-harmonic polynomials of degree $m$ of the form
$$e^{in\theta}\rho^m\sin^{m/2}\phi h(cot\phi)$$
where $|n|\leq m,\;m\equiv n$ (mod 2) and $h(x)$ is some polynomial in $x$  which, depends on $\alpha,\;\; n$ and $m$.\\
\\
\textbf{Proof:}~
We note that
$${\cal L}_{\alpha}\left(e^{in\theta}f(|z|,t)\right)=e^{in\theta}{\cal L}_{\alpha,n}f(|z|,t),\eqno (2.8)$$
where ${\cal L}_{\alpha,n}$ is a second order differential operator in $r=|z|$ and $t$.  Let $u(z,t)$ be a $\h_1$-homogeneous ${\cal L}_{\alpha}$-harmonic polynomial of degree $m$.  We may write it, uniquely as follows.
$$u(z,t)=\sum_ne^{in\theta}u_n(|z|,t).$$
Then
$${\cal L}_{\alpha}(u)=\sum_ne^{in\theta}{\cal L}_{\alpha,n}\;u_n(|z|,t)=0, \eqno (2.9)$$
which implies that ${\cal L}_{\alpha,n}u_n(|z|,t)=0$ for all $n$.  Therefore, from $(2.8),$ we have
$${\cal L}_{\alpha}(e^{in\theta}u_n(|z|,t))=0, $$
for all $n$.  The general term of $u$ is some constant times
$$z^{m_1}\;\bar{z}^{m_2}\;t^{m_3}=e^{i(m_1-m_2)\theta}|z|^{m_1+m_2}t^{m_3}=e^{in\theta}|z|^m
\left(\frac{t}{|z|^2}\right)^{m_3}$$
where $m_1+m_2+2m_3=m$ and $n=m_1-m_2$.  In particular $m\equiv n(mod\; 2),$\\
$\;|n|\leq m$, and since $|z|^m=\rho^m\sin^{m/2}\phi$ and $t/|z|^2=\cot\phi$, the proof is complete. \hfill $\Box$\\
\\
\textbf{Remark 2.2.3}~~In the calculation of $h(x)=h^{^{(\alpha,n)}}_{_{(m-|n|)/2}}(x),$ an explicit computation of
$${\cal L}_{\alpha}(u)=-Z\bar{Z}u+i(\alpha-1)\frac{\partial u}{\partial t}$$
is required,
where
$$u=e^{in\theta}\rho^m\sin^{m/2}\phi y(cot\phi).$$
A computation in [9] shows that
$$\begin{array}{rcl}
{\cal L}_{\alpha}(u)&=&-e^{in\theta}\rho^m\sin^{m/2}\phi\{(1+cot^2\phi)y^{\prime\prime}(cot\phi)\\
&&-[(m-1)cot\phi+i(\alpha+n)]y'(cot\phi)\\
&&+\frac{1}{4}(m^2-n^2)y(cot\phi)\}\end{array}$$
We assume ${\cal L}_{\alpha}(u)=0$ and then look for a polynomial solution, $y(x)$ of
$$(1+x^2)y^{\prime\prime}(x)-(i(\alpha+n)+(m-1)x)y'(x)+\frac{1}{4}(m^2-n^2)y(x)=0. \eqno (2.10)$$
To find the required solution, we substitute
$$y=\sum^{\infty}_{v=0}a_v(x-i)^v\eqno (2.11)$$
in (2.10).  This leads to
$$(x+i)\;\sum^{\infty}_{v=2}v(v-1)a_v(x-i)^{v-1}(i(\alpha+n)+(m-1)x\sum^{\infty}_{v=1}v a_v(x-i)^{v-1}$$
$$+\;\frac{m^2-n^2}{4}\sum^{\infty}_{v=0}v a_v(x-i)^v=0,$$
which yields the recurrence relation
$$i(2v-\alpha-n-m+1)(v+1)a_{v+1}+(v^2-mv+\frac{m^2-n^2}{4})a_v=0\eqno (2.12)$$
for $v=0,1,2,3,\cdots$.  Now $m+n$ is even and by hypothesis $\pm\alpha=1,3,5,7,\cdots$.  Therefore, the coefficient of $a_{v+1}$ in $(2.12)$ never vanishes.  Thus $y(x)$ is a polynomial if and only if the coefficient of $a_v$ in $(2.12)$ vanishes at some point, i.e.,
$$0=v^2-mv+\frac{m^2-n^2}{4}=\left(v-\frac{m-|n|}{2}\right)\left(v-\frac{m+|n|}{2}\right).$$
Since $|n|\leq m$, $(2.10)$ has a polynomial solution, which is unique up to a constant multiplier, of degree
$k=\frac{1}{2}(m-|n|),$ we set $a_0=1$ in $(2.11).$ Then
$$a_v=\left(\frac{i}{2}\right)^v\;\frac{\left(-\frac{m-|n|}{2}\right)_v\left(-\frac{m+|n|}{2}\right)_v}
{\disp\left(\frac{-m+n}{2}-\frac{\alpha-1}{2}\right)_v}\;\cdot\;\frac{1}{v!},$$
$v=1,2,3,\cdots$ where we used the standard notation
$$a_0=1, (a)_n=a(a+1)\cdots(a+n-1),\;n=1,2,3,\cdots.$$
Thus the polynomial solution of (2.10) is given by
$$\begin{array}{rcl}
y(x)&=&\disp\sum^{\frac{1}{2}(m-|n|)}_{v=0}\;\left(\frac{i}{2}\right)^v\;\frac{\left(\frac{-m-|n|}{2}\right)_v
\;\left(-\frac{m+|n|}{2}\right)_v}{\left(-\frac{m+n}{2}-\frac{\alpha-1}{2}\right)_v}\\
\ \\
&=&F\left(-\frac{m-|n|}{2},\;\frac{m+|n|}{2}; - \frac{m+n}{2}-\frac{\alpha-1}{2};\frac{1}{2}+\frac{1}{2}ix\right)
.\end{array}\eqno (2.13)$$
Here $F(a,b;c;z)$ stands for Gauss' hypergeometric function defined by
$$F(a,b; c;z) = \sum^{\infty}_{v=0}\frac{(a)_v(b)_v}{(c)_v}\;\frac{z^v}{v!}$$
as long as $c\neq 0, -1, -2, -3,\cdots$.  We note that in our case
$$c = - \frac{m+n}{2}-\frac{\alpha-1}{2}\neq0,-1,-2,-3,\cdots$$
since $\pm\alpha= 1, 3, 5,\cdots$.  For sufficiently small $z$,
$$F(a,b; c;z) = (1-z)^{-a}\;F\left(a,c-b; c; \frac{z}{z-1}\right)\;,\;\eqno (2.14)$$
as long as $c\neq 0, -1, -2, -3,\cdots,$ [11(I)].  If $a=0, -1, -2, -3,\cdots$ both sides of $(2.14)$ are polynomials, hence $(2.14)$ holds for all $z\in \C$.  In our case, $a=-k=-(m-|n|/2)$ is a non-positive integer, therefore, disregarding constant multipliers, $(2.13)$ can be written in the form
$$y(x)=(x+i)^{(m-|n|)/2}F\left(-\frac{m-|n|}{2},\;\frac{|n|-n}{2}-\frac{\alpha-1}{2};\frac{m+n}{2} - \frac{\alpha-1}{2};\frac{x-i}{x+i}\right).$$
This formula can be further simplified in the following way.  With a non-negative integer, we have\\
$\disp\left(\begin{array}{c}
k+b\\
k\end{array}\right)(x+i)^k\;F\left(-k, -k-a; b+1; \frac{x-i}{x+i}\right)$\\
$$\begin{array}{rcl}
&=&\left(\begin{array}{c}
k+b\\
k\end{array}\right)(x+i)^k\disp\sum^k_{v=0}\frac{k(k-1)\cdots(k-v+1)}{(b+1)(b+2)\cdots(b+v)}
\left(\begin{array}{c}
k+a\\
v\end{array}\right)\left(\frac{x-i}{x+i}\right)^v\end{array}$$
$$\begin{array}{rcl}
&=&\disp(x+i)^k\sum^k_{v=0}\frac{(k+b)(k+b-1)\cdots(k+v+1)}{1.2\cdots(k-v)}
\left(\begin{array}{c}
k+a\\
v\end{array}\right)\left(\frac{x-i}{x+i}\right)^v\\
\ \\
&=&\disp\sum^k_{v=0}\left(\begin{array}{c}
k+a\\
v\end{array}\right)\left(\begin{array}{c}
k+b\\
k-v\end{array}\right)(x-i)^v(x+i)^{k-v}.\end{array}$$
Thus, with $k=(m-|n|)/2, k+a=-(|n|-n)/2+(\alpha-1)/2$ and\\ $b=-(m+n)/2-(\alpha+2)/2$, we the following result.\
\\
\textbf{Proposition 2.2.4}~
$h(x)=h^{^{(\alpha,n)}}_{_{(m-|n|)/2}}(x)$ with $(m-|n|)/2=0,1,2,\cdots$ is the unique, up to constant multiplier
, polynomial solution of (3.25).  $h^{^{\alpha,n)}}_{_{k}}(x), \; k=0,1,2,\cdots$ is exactly of degree $k$.  Explicitly, one has
$$h^{^{(\alpha,n)}}_{{k}}(x)=\sum^k_{v=0}\left(\begin{array}{c}
-\frac{|n|+n}{2}-\frac{\alpha+1}{2}\\
\ \\
v\end{array}\right)\;\left(\begin{array}{c}
-\frac{|n|-n}{2}+\frac{\alpha-1}{2}\\
\ \\
k-v\end{array}\right)\;(x+i)^v(x-i)^{k-v}.$$\\
\\
\textbf{Remark 2.2.5}

Using propositions $2.2.2$ and $2.2.4,$ the basis elements of the solid ${\cal L}_{\alpha}$-spherical harmonics of degree $m$ have the form
$$e^{in\theta}\rho^m\sin^{m/2}\phi h^{^{(\alpha,n)}}_{_{(m-|n|)/2}}(\cot \phi),$$
where $n=0,\pm 1,\pm 2, \pm 3,\cdots$ such that $(m-|n|)/2=0,1,2,3,\cdots$.

Now, to simplify this, let $x=cot \phi.$ Then
$$x+i=\frac{\cos\phi}{\sin\phi}+i=\frac{e^{i\phi}}{\sin\phi},\;\; x-i=\frac{e^{i\phi}}{\sin\phi}.$$
Therefore,
$$\sin^k\phi h^{^{(\alpha,n)}}_{{k}}(cot \phi)=\sum^k_{v=0}\left(\begin{array}{c}
-\frac{|n|+n}{2}-\frac{\alpha+1}{2}\\
\ \\
v\end{array}\right)\;\left(\begin{array}{c}
-\frac{|n|-n}{2}+\frac{\alpha-1}{2}\\
\ \\
k-v\end{array}\right)\;e^{iv\phi}e^{-i(k-v)\phi}$$
$$=\sum^k_{v=0}\left(\begin{array}{c}
-\frac{|n|+n}{2}-\frac{\alpha+1}{2}\\
\ \\
v\end{array}\right)\;\left(\begin{array}{c}
-\frac{|n|-n}{2}+\frac{\alpha-1}{2}\\
\ \\
k-v\end{array}\right)\;e^{i(2v-k)\phi}.$$
We set
$$H^{^{(\alpha,n)}}_{_{k}}(e^{i\phi})=(-1)^k\sum^k_{v=0}\left(\begin{array}{c}
-\frac{|n|+n}{2}-\frac{\alpha+1}{2}\\
\ \\
v\end{array}\right)\;\left(\begin{array}{c}
-\frac{|n|-n}{2}+\frac{\alpha-1}{2}\\
\ \\
k-v\end{array}\right)\;e^{i(2v-k)}.$$
Then,
$$\sin^k\phi h^{(\alpha,n)}_k(cot\phi) = (-1)^kH^{(\alpha,n)}_k(e^{i\phi}).$$
We now state the main theorem of this section.\\
\\
\textbf{Theorem 2.2.6}~
Set
$${\cal L}_{\alpha}=-Z\bar{Z}+i(\alpha-1)\frac{\partial}{\partial t},\;\;\pm\alpha=1,3,5,\cdots .$$
Then the linear space ${\cal H}_M^{(\alpha)}$ of $\h_1$-homogeneous ${\cal L}_{\alpha}$-harmonic polynomials of degree $m,\;\; m=0,1,2,\cdots$ has dimension $m+1$.  A basis for ${\cal H}_m^{(\alpha)}$ can be found as follows:  For each $n$ and $k$, such that $m=2k+|n|$ with $k=0,1,2,3,\cdots$ and $n=0,\pm 1,\pm 2, \pm 3,\cdots$ set
$${\cal H}_m^{(\alpha,n)}(\rho,\theta,\phi)=e^{in\phi}\rho^m\sin^{|n|/2}\phi\;H^{^{(\alpha,n)}}_{_{(m-|n|)/2}}
(e^{i\phi}),.$$
where $\disp H^{^{(\alpha,n)}}_{_{k}}(e^{i\phi})$ is defined by the generating function
$$(1-\rho e^{i\phi})^{-(|n|+n)/2-(\alpha+1)/2}(1-\rho e^{-i\phi})^{-(|n|-n)/2+(\alpha-1)/2}$$
$$=\sum^{\infty}_{k=0}\rho^kH^{(\alpha,n)}_k(e^{i\phi}),$$
with the right hand side converging for $\rho<1$.

Then $\left\{{\cal H}^{^{(\alpha,n)}}_{_{m}}\right.$, all possible $\left.n\right\}$ is a basis for ${\cal H}^{^{(\alpha)}}_{_{m}}$.\\
\\
\textbf{Definition 2.2.6}~
The ${\cal L}_{\alpha}$-spherical harmonics are obtained by restricting $\h_1$-homogeneous ${\cal L}_{\alpha}$-harmonic polynomials to $bB_{\h_1}(1)$, i.e., to $\rho=1$.\\
\\
\textbf{Theorem 2.2.7}~  Each ${\cal L}_{\alpha}$-spherical harmonic is a unique linear combination of functions of the form
$$e^{in\theta}\sin^{|n|/2}\phi H^{(\alpha,n)}_k(e^{i\phi})$$
 with $k=0,1,2,3,\cdots,$ and $n=0,\pm1, \pm 2,\pm 3,\cdots.$\\
\textbf{ Proof:}~
We recall the homogeneous harmonic polynomials of $m,\; m=0,1,2,\cdots$ of the Laplacian in $\R^3$.  In spherical coordinates they are
$$e^{in\theta}\;r^m\sin^{|n|}\phi\;P^{^{|n|+1/2}}_{_{m-|n|}}(\cos\phi),\eqno (2.15)$$
$0\leq \theta\leq 2\pi,\;0\leq\phi\leq\pi$ and $0<r<\infty$, where $n=0,\pm1,\pm2,\cdots,\pm m$ and the
ultraspherical (or Gegenbauer) polynomial, $P^2_k(x), -1\leq x \leq$, are defined via the generating function
$$\begin{array}{rcl}
(1-re^{i\phi})^{-\lambda}(1-re^{-i\phi})^{-\lambda}&=&(1-2r\cos\phi+r^2)^{-\lambda}\\
\ \\
&=&\disp\sum^{\infty}_{k=0}r^kP^{\lambda}_k(\cos\phi).\hspace{2.0in}\end{array} (2.16)$$
Sometimes $P^{\lambda}_k$ are denoted y $C^{\lambda}_k$ [10].  The central result in the theory of the polynomials $P^{\lambda}_k(x)$ is\\
\\
\textbf{Proposition 2.2.8}~
For each $\lambda >-\frac{1}{2}$ the polynomials $P^{\lambda}_k(x), \;k=0,1,2,3,\cdots$ form an orthogonal basis for $L^2([-1,1])$, with measure $(1-x^2)^{\lambda-\frac{1}{2}}dx$.  Furthermore,
$$\begin{array}{rcl}
\disp\int^1_{-1}(P^{\lambda}_k(x))^2(1-x^2)^{\lambda-\frac{1}{2}}dx &=& \frac{\pi^{\frac{1}{2}}(2\lambda)_k\Gamma(\lambda+\frac{1}{2})}{(k+\lambda)k!\Gamma(\lambda)},\\
\ \\
k&=&0,1,2,3,\cdots\end{array}\eqno (2.17)$$
(For proof, see [11(II)].
It follows from C(17) that with $\ell=0,1,2,3,\cdots,$
$$P^{^{\ell+\frac{1}{2}}}_{_{k}}(\cos\phi) = H^{^{(-2\ell,-2\ell)}}_{_{k}}(e^{i\phi}),\;\;k=0,1,2,\cdots.$$
In particular, when $\ell=0$, we obtain the Legendre polynomials
$$P_{_{k}}(\cos\phi)=P^{^{\frac{1}{2}}}_{_{k}}(\cos\phi)=H^{^{(0,0)}}_{_{k}}(e^{i\phi}),\;\;k=0,1,2,
\cdots.\eqno (2.18)$$

\begin{center}{\bf References}\end{center}
\begin{description}

\item[] [1]~Andrews, G.E., Askey, R. and  Roy, R. 2000.  Special functions. \emph{Encyclopaedia of Math. and its Appl., Cambridge Univ. Cambridge}.
\item[][2]~	Howe R.: on the role of the Heisenberg group in harmonic analysis.  \emph{Bull. Amer. Math. Soc.,} 3(1980), 821-843.
\item[][3]~  Stein E.M.: Harmonic analysis: Real variable methods, orthogonality and Oscillatory Integrals. \emph{Princeton Univ. Press Princeton,} 1993.
\item[][4]	~Folland G.B.: A fundamental solution for a subelliptic operator. \emph{Bull. Amer. Math. Soc.,} 79(2) (1973), 373-376. London Math.Soc. Lecture Note Series, Cambridge Univ. Press (1979).
\item[][5]~	Rothschild L.P.: Local solvability of left-invariant differential operators on the Heisenberg group. \emph{Proc. Amer. Math. Soc.,} 74(2) (1979), 383-388.
\item[][6]~ Bell, W.W. 1968. Special functions for Scientist and Engineers. \emph{D. Van Nostrand Co. Ltd., Princeton.}
\item[][7]~    Whittaker, E.T. and Watson, G.N.:  A Course of modern analysis.\emph{ Cambridge Univ. Press Cambridge,  1996.}
\item[][8] ~Folland, G.B. and Stein, E.M.:  Estimate for the   complex and analysis on the Heisenberg group. \emph{Comm. Pure and Appl. Math.,} XXVII(1974) , 429-522.
\item[][9]  ~Greiner, P.C. 1980.  Spherical Harmonics on the Heisenberg group. \emph{Canad. Math. Bull., 23 (4) (1980), 383-396.}

\item[][10]~ Egwe M.E.: Aspects of Harmonic analysis on the Heisenberg group. Ph.D. thesis, University of Ibadan, Ibadan, Nigeria, 2010.

\item[][11]~ Egwe M.E.: Non-Solvability of Heisenberg Laplacian by Factorization, \emph{Journal of Mathematical Sciences, 21 (1)(2010), 11-15.}

\item[][12]~ Egwe M.E. and U.N. Bassey: On Isomorphism Between Certain Group Algebras on the Heisenberg Group, \emph{J. Math. Phy. Anal. Geom. 9(2) (2013), 150-164.}

\item[][13]~     Erdelyi, A.:  Higher transcendental functions, Vols. I \& II. \emph{Robert E. Krieger Publ., Florida,  1953.}.
\end{description}

\end{document}